\documentclass[a4paper]{easychair}
\usepackage[utf8]{inputenc}
\usepackage{rotating}
\usepackage{pdflscape}
\newcommand{\equivalence}[2]{{#1} \leftrightarrow {#2}}
\newcommand{\implicationsymbol}[0]{\rightarrow}
\newcommand{\implication}[2]{{#1} \implicationsymbol {#2}}
\newcommand{\conjunctionsymbol}[0]{\wedge}
\newcommand{\conjunction}[2]{{#1} \conjunctionsymbol {#2}}
\newcommand{\disjunctionsymbol}[0]{\vee}

\newcommand{\bottom}[0]{\perp}

\newcommand{\consequencesymbol}[0]{\models}
\newcommand{\consequence}[2]{{#1} \consequencesymbol {#2}}
\let\phi=\varphi
\usepackage[utf8]{inputenc}
\usepackage{amsmath}
\usepackage{amssymb}

\newcommand{\provabilitysymbol}[0]{\mathtt{P}}
\newcommand{\negationsymbol}[0]{\neg}
\newcommand{\negate}[1]{\negationsymbol{#1}}
\newcommand{\see}[1]{\mathsf{C}_{#1}}
\newcommand{\seeone}[0]{\see{1}}
\let\phi=\varphi

\usepackage{hyperref}
\hypersetup{pdftitle=A machine-assisted view of paraconsistency}
\hypersetup{pdfauthor=Jesse Alama}
\hypersetup{pdfsubject={paraconsistency}}
\hypersetup{pdfkeywords={paraconsistent logics, automated theorem proving, theory exploration}}
\begin{document}
\title{A machine-assisted view of paraconsistency}
\titlerunning{Machine-assisted paraconsistency}
\author{Jesse Alama\inst{1}}
\institute{Theory and Logic Group{\\}Technical University of Vienna{\\}\email{alama@logic.at}}
\authorrunning{Alama}

\clearpage

\maketitle

\begin{abstract}{For a newcomer, paraconsistent logics can be difficult to grasp.  Even experts in logic can find the concept of paraconsistency to be suspicious or misguided, if not actually wrong.  The problem is that although they usually have much in common with more familiar logics (such as intuitionistic or classical logic), paraconsistent logics necessarily disagree in other parts of the logical terrain which one might have thought were not up for debate.  Thus, one's logical intuitions may need to be recalibrated to work skillfully with paraconsistency.  To get started, one should clearly appreciate the \emph{possibility} of paraconsistent logics and the genuineness of the distinctions to which paraconsistency points.  For this purpose, one typically encounters matrices involving more than two truth values to characterize suitable consequence relations.  In the eyes of a two-valued skeptic, such an approach might seem dubious.  Even a non-skeptic might wonder if there's another way.  To this end, to explore the basic notions of paraconsistent logic with the assistance of automated reasoning techniques.  Such an approach has merit because by delegating some of the logical work to a machine, one's logical ``biases'' become externalized.  The result is a new way to appreciate that the distinctions to which paraconsistent logic points are indeed genuine.  Our approach can even suggest new questions and problems for the paraconsistent logic community.}\end{abstract}

\section{Opening non-paraconsistent minds}\label{sec:introduction}

Given a set $X$ of sentences (formulas, statements, propositions), a notion $\models$ of logical consequence (a relation between subsets of $X$ and $X$) is said to be \emph{explosive} if it satisfies the condition:
\[
\text{$\consequence{\{\, \phi, \negate{\phi} \,\}}{\psi}$ for all $\phi$ and $\psi$ in $X$}
\]

{\noindent}The consequence relation $\models$ is said to be \emph{paraconsistent} if it is not explosive.

Explosive consequence relations are easily found.  The standard consequence relations for classical logic, intuitionistic logic, and all intermediate logics are explosive.  Even ``alternative'' logics such as modal logics are explosive.  Everyone knows that from explicitly contradictory information one can infer whatever one wants.  Right?

From this perspective, paraconsistency seems impossible, so the claim that it \emph{is} possible appears unorthodox, contrarian, or even heretical.  It seems that the meaning of certain basic logical terms such as ``negation'', ``contradiction'', or ``consistency'' has to be changed if the whole approach is to even make sense~\cite{slater1995paraconsistent}.

This is only the expression of an attitude, or perhaps just a complaint; it is not an impossibility proof.  By keeping an open mind and not rejecting paraconsistency out of hand, one might ask a mild question: \emph{Under what conditions is paraconsistency possible?}  Various philosophical arguments have been put forward to justify the sensibility of paraconsistency~\cite{brown1999yes,tanaka2013making}.  They point to apparent examples of paraconsistency in some of our reasoning.  Even adopting a rigorous approach, one can point to various paraconsistent logics (e.g., da~Costa's hierarchy~$\see{k}$~\cite{alves1977semantical})  So paraconsistency is, evidently, possible, at least from the perspective of abstract consequence relations.  But perhaps these philosophical arguments for the plausibility of paraconsistent, and semantic approaches demonstrating the possibility of paraconsistent logics, are somehow unsatisfying.  Is there another way?

One avenue for taking a fresh approach is to treat paraconsistency as a target for automated theory exploration.  More precisely, paraconsistency becomes a problem for knowledge representation, to be investigated with automated theorem provers (ATPs) and model (or countermodel) finders .  On this approach, one ``externalizes'' one's assumptions about logic by stating them explicitly in a logical language, and then uses automated reasoning systems to explore the subject by trying to prove or refute statements about paraconsistent logic that have been rigorously formulated.   (To be clear, we mean to use ATPs to reason \emph{about} paraconsistency, not to reason \emph{within} a paraconsistent logic~\cite{dacosta1990automatic}.)

The approach is hardly new; see, for example, the work of Ulrich~\cite{ulrich2001legacy} and Wos~\cite{wos2003automated} for many similar applications, e.g., to intuitionistic logic, classical logic, relevance logics, modal logics, etc.  To our knowledge, though, the idea of exploring the foundations of paraconsistency with the help of automated theorem provers, in the spirit of Lakatos, is new.  We believe that such an approach can complement standard semantic methods for investigating paraconsistency, offers a fresh argument for the possibility of paraconsistency, and can even suggest new problems in research into paraconsistency.

After considering the impact our choice of representation may have on the interpretation of our results (Section~\ref{sec:representation}), the heart of the paper consists of a series of (machine-assisted) experiments bearing on paraconsistency are considered.  In Section~\ref{sec:experiment-0} we consider the very possibility of paraconsistency.  We work with da~Costa's~$\seeone$~\cite{dacosta1974theory}.  One starts to glimpse what paraconsistency has to offer by finding that the one cannot prove in~$\seeone$ that, for every $\phi$, that $\conjunction{\phi}{\negate{\phi}}$ is explosive.  In Section~\ref{sec:experiment-1} we consider the amusing case of trivializing triviality.  In Section~\ref{sec:experiment-3} we take an opposite tack, focusing on the conditions for the possibility of explosiveness (or non-paraconsistency).  Finally, in Section~\ref{sec:experiment-4} we concern ourselves with the existence of $\bottom$-like explosive particles.  Section~\ref{sec:conclusion} summarizes and considers possibility for further work along the lines taken here.

\section{Representation: first-order provability structures}\label{sec:representation}

Our target is a class of paraconsistent propositional logics.  We focus, more specifically, on da~Costa's~$\seeone$~\cite{dacosta1974theory,carnielli2007logics}.  To get started, we postulating a primitive unary predicate $\provabilitysymbol$, intended to represent the set of theorems.  The intended domain of discourse consists of propositions.  Our aim is to lay down axioms that capture a set of theorems.

It may seem hopeless to take such an abstract approach, but such an approach pays off.  By adopting an approach that is amendable to ATPs, we put ourselves in a position to reap their fruit.  Notwithstanding some important expressive limitations of first-order logic, ATPs have the delightful ability to, one the one hand, confirm that what we expect to be provable is in fact provable, but perhaps in a way which is wholly unexpected.  On the other hand, an ATP is often able to reveal genuine gaps in our own thinking.  If one is new to paraconsistency, one might even find that an ATP can assist one to see the field with an open mind.

To represent consequence relations (relations $\models$ between subsets of a set $X$ and $X$), we work with elementary classes of first-order structures in a language with a distinguished unary predicate~\verb+P+, intended to represent the property of being provable from no assumptions (that is, the property of being a logical truth in whatever logic we are axiomatizing).  The intended domain of discourse consists of propositions.  The languages we consider include some connectives, such as implication~\verb+i+ and negation~\verb+n+.  One might also consider adding disjunction~\verb+o+ (``or'') and conjunction~\verb+a+ (``and'').  Using~\verb+P+ and these connectives, one can axiomatize various logics by postulating various formulas as theorems.

It is important to reflect on the properties that consequence relations have in virtue of the representation we have chosen.


\subsection{Uniform substitution}

In the study of paraconsistency, one may need to pay close attention to a feature of most propositional logics that is often tacitly used many times.  Namely, the idea that, from the standpoint of provability, atoms can be thought of as variables into which one can substitute values (other formulas) while preserving provability.  To take the simplest case, if an atom $p$ occurs in a theorem $\phi$, then, for any formula $\psi$, the instance $\phi[p \coloneqq \psi]$ of $\phi$ is also a theorem.  In our representation, we postulate modus ponens; but owing to the way that many automated theorem provers work, we are effect postulating condensed detachment, which is a fusion of modus ponens and unification into a single rule of inference~\cite{kalman1983condensed}.  Substitution is thus ``built-in''.  See also McCune and Wos~\cite{mccune1992experiments}.

In the study of paraconsistency, for various reasons uniform substitution may need to be abandoned, though there are also paraconsistent logics for which uniform substitution does hold.  In our representation, uniform substitution does hold.

\subsection{Identity of propositions}

Given an atomic formula $p$, the implication $\implication{p}{p}$ is, obviously, not the same thing as $p$.  From a logical point of view, the former is contingent whereas the second is a logical truth; from a syntactic point of view, the former has length $1$, whereas the latter has length $3$ (counting all symbols).  There are, moreover, infinitely many distinct formulas: $p$ is not the same as $\implication{p}{p}$, which is not the same as $\implication{p}{\implication{p}{p}}$, etc. Of course, some formulas turn out to be equivalent to one another and so perhaps could be identified.  But at least from a syntactic point of view, there are, obviously, infinitely many distinct formulas.

When logics are represented in the first-order manner that we have adopted in this paper, one soon encounters the perplexing fact that there exist finite models of our axioms.  Since the intended domain of discourse consists of formulas, and since implication and negation are present, how can it be there are only finitely many formulas?  Why, in other words, can't we replay the argument in the preceding paragraph to the effect that there are infinitely many formulas?

From a first-order model of~$\seeone$ (especially, finite ones), one can construct a structure which perhaps accords more one's intuition.  What is wanted is a structure $(X, \vdash)$ where $X$ is an infinite set of formulas and $\vdash$ is a distinguished subset of $X$ (the theorems).  $X$ should, moreover, be \emph{generated} by formula-building functions (connectives) such as $\implicationsymbol$, $\negationsymbol$, $\conjunctionsymbol$, $\disjunctionsymbol$.  The difficulty is to construct a suitable homomorphism from an arbitrary algebra of formulas into a an \emph{generated} algebra of formulas.  We do not go into the details here.

\subsection{Expressive limitations of FOL}

Using the plain language of FOL in which we represent our problems, a few natural properties of structures are in fact not representable.  We are unable to express, for example, that the domain of discourse (formulas) are \emph{generated} by a set of atoms.  Only partial remedies are avaiable.  We can define the notion of atom as not being the value of the negation, conjunction, disjunction, or implication formula-building operations and then postulate, for instance, the existence of atoms.  We can even define the notion of being an implication, being a negation, etc., and postulate that everything is either (1) an atom, or (2) an implication, or (3) a conjunction, or\dots.  But such efforts still fail to capture the familiar property of propositional languages that they are generated by their atoms, that negation is an injection into the set of formulas, etc.

We can postulate principles that ensure that the domain of discourse is infinite.  For exbample, dealing with an implication and negation, we can ensure that the domain of discourse is infinite by postulating that

\begin{itemize}
\item Implications are not negations;
\item The negation of $\phi$ is never $\phi$ (that is, $\negationsymbol$ is an injection on the set of formulas);
\end{itemize}

{\noindent}Given an implication $\phi$, these conditions ensure that $\negate{\phi}$, $\negate{\negate{\phi}}$, \dots, are all distinct objects.  Interestingly, for the purposes of proving theorems about paraconsistent logics these assumptions are usually not needed and can be safely dropped.  From the standpoint of finite model-building this fact is evidently advantageous because it makes clear the possibility of working with finite models.

A further consequence of the expressive weakness of first-order logic concerns the generation of theorems (as opposed to the generation of formulas).  The familiar distinction between \emph{logic} and \emph{theory in a logic} is that a logic lays down the groundwork for consequence, whereas a theory postulates additional axioms while still working within the logic.  From a mathematical point of view, a (propositional) logic consists of a set of formulas generated by a set of atoms and some connectives, and closed under a set of inference rules.  A consequence of expressive limitations of FOL is that we generally cannot give a first-order formula that captures the fact that the set of theorems is generated.  So when examining various models, we find that there can be ``theorems from nowhere'', that is, we can be working with a theory rather than just a logic.


\subsection{The deduction theorem}

The notion of provability-from-a-single-assumption is, in our case, defined in terms of the notion of provability and implication:

\[
\phi \models \psi \quad \text{iff} \quad \models \implication{\phi}{\psi} .
\]

Whereas for our purposes this is a definition, normally this property of provability is has to be proved (and is known as the deduction theorem).  We thereby exclude the possibility of paraconsistent logics in which the deduction theorem does not hold.  Similarly, the notion of provability-from-two-assumptions is defined in terms of provability-from-one-assumption in such a way as to preserve the deduction theorem:

\[
\phi,\psi \models \chi \quad \text{iff} \quad \phi \models \implication{\psi}{\chi} .
\]

{\noindent}There is apparently some non-uniqueness here, but that is addressed by the fact that~$\seeone$ proves

\[
\models \implication{\phi}{(\implication{\psi}{\chi})} \quad \text{iff} \quad \models \implication{\psi}{(\implication{\phi}{\chi})}.
\]




\section{Experiment 0: The very possibility of paraconsistency}\label{sec:experiment-0}

Working with da~Costa's~$\seeone$~\cite{dacosta1974theory}---that is, representing~$\seeone$ in terms of the undefined predicate $\provabilitysymbol$ and a few sentence-building functions for implication, conjunction, disjunction, and negation---one starts to glimpse what paraconsistency has to offer by finding that the one cannot prove in~$\seeone$ that, for every $\phi$, that $\conjunction{\phi}{\negate{\phi}}$ is explosive.  Very simple models illustrating this possibility are available.

\section{Experiment 1: Trivializing triviality}\label{sec:experiment-1}

Call a proposition \emph{conditionally explosive} if its provability entails the provability of all propositions, that is, it its provability trivializes the logic.  The definition is, in other words,

\[
\text{If $\models{\phi}$, then for all $\psi$ we have $\models{\psi}$}
\]

Are there conditionally explosive propositions?  The answer is, briefly, yes.  With the assistance of a machine an unexpectedly clever answer (namely, Answer~$4$ below) to the question can be found.

\paragraph{Answer 1.} Is there a $\bottom$?  If there were a distinguished formula $\bottom$ with a standard $\implication{\bottom}{\phi}$, then of course we would be finished.  It is consistent with~$\seeone$ that there is a $\bottom$, but we have not postulated such a formula, so we have to keep working.

\paragraph{Answer 2.} Even without a $\bottom$, we can perhaps answer the question with only marginally more work.  Take a formula $\phi$ and consider $\conjunction{\phi}{\negate{\phi}}$.  This may not work, though, as the results in Section~\ref{sec:experiment-0} suggest.  In many logics (e.g., intuitionistic logic, all intermediate logics, classical logic, etc.) such a formula answers the question.  It is consistent with~$\seeone$ that there is such a formula, but the results of Section~\ref{sec:experiment-0} tell us that there's no hope of \emph{proving} the existence of such a formula.

\paragraph{Answer 3.} Uniform substitution does the job.  Consider an atomic formulas $p$.  If $p$ were provable, then by uniform substitution all formulas $\phi$ would be provable.

\paragraph{Answer 4.} For any unary predicate symbol $P$, the following is a (somewhat awkward) theorem of classical first-order logic:

\[
\exists x [ \implication{P(x)}{\forall y P(y)} ]
\]

{\noindent}(To see that it is provable, consider its negation.) Now consider the consequence relation $\models$ of~$\seeone$ and the definition of the predicate $e(\phi)$ (``$\phi$ is explosive''):

\[
\forall x [ \equivalence{e(x)}{(\implication{\models x}{\forall y \models y})} ]
\]

The question is to prove or refute the following sentence:

\[
\exists x [ \implication{\models x}{\forall y \models y} ]
\]

This formula is indeed provable: it is just a recapitulation (or rather, an instance) of the curious theorem above

What is amusing here is that no principles of $\models$ were used.  Thus, we needn't postulate that~$\seeone$ is closed under modus ponens, or, really, anything at all about it.  The definition of conditionally explosive formula alone does \emph{all} the work.  We have proved the existence of conditionally explosive formulas from, in effect, nothing.

The result of this experiment is that we need to distinguish conditionally explosive propositions from $\bottom$-like formulas; see Section~\ref{sec:experiment-4}.

\section{Experiment 2: The possibility of explosiveness}\label{sec:experiment-3}

The logic~$\seeone$ is a basis for paraconsistency, as confirmed in Section~\ref{sec:experiment-0}.  But~$\seeone$ is not essentially paraconsistent in the sense that it is consistent with~$\seeone$ that the explosion principle

\[
\implication{p}{(\implication{\negate{p}}{q})}
\]

{\noindent}is provable.  If the explosion principle is postulated, then one can simply delete the axioms of~$\seeone$ governing the consistency operator $^{\circ}$ because they all become dependent.  Interestingly, even if the explosion principle is postulated, one cannot delete the principle of excluded middle, thereby reconfirming the idea that explosion is not inherently tied to classical logic.

\section{Experiment 3: Existence of $\bottom$-like formulas}\label{sec:experiment-4}

Call a formula $\phi$ \emph{$\bottom$-like} if it has the property:
\[
\text{$\models {\implication{\phi}{\psi}}$ for all formulas $\psi$}
\]

Do there exist $\bottom$-like formulas?  One can verify easily that their existence is consistent with~$\seeone$.  Further, one can postulate the existence of $\bottom$-like formulas and still remain paraconsistent.

Curiously, it is somewhat difficult to \emph{avoid} $\bottom$-like formulas.  At the time of writing, the author has so far found no $\bottom$-free~$\seeone$-theories have uncovered by our search methods.  The difficulty is that small finite models of~$\seeone$ all have at least one bottom-like formula.  It is possible that every finite model of~$\seeone$ does have a bottom-like formula and there exists an infinite model of~$\seeone$ without a bottom-like particle.

\section{Conclusion and future work}\label{sec:conclusion}

We have considered the concept of paraconsistency from the perspective of a machine, trying to externalize our assumptions about logic and investigating paraconsistency systematically with the help of automated reasoning technology.  We have motivated the basic distinctions in paraconsistency, and even offered a fresh argument for its possibility, by pursuing this route.

Like da~Costa, the consistency operator $^{\circ}$ is in our approach a defined operation on formulas.  It could, however, be considered as primitive.  Apropos Section~\ref{sec:experiment-0}, one could even consider the possibility of other definitions of $^{\circ}$ that accomplish precisely what it does.  In the context of so-called lightly explosive logics, one might consider the following principle (which one might call ``imminent explosions''):

\[
\text{For all $\phi$, there exists a $\psi$ such that for all $\chi$, we have $\{ \phi, \negate{\phi}, \psi \} \models \chi$}
\]

{\noindent}which is an abstract variation of the notion of finite gentle explosion~\cite{carnielli2007logics}.

Apropos Section~\ref{sec:experiment-3}, it would be worth investigating whether intuitionistic logic is obtained by postulating explosion but removing the principle of excluded middle.  In Section~\ref{sec:experiment-4} we were unable to settle the question of whether $\bottom$-like particles always exist.  Presumably, they do not always exist (working, that is, in~$\seeone$), but so far our methods have not been able to resolve this one way or another.

In the end, one might take the view that paraconsistency is simply impossible or incoherent.  The most to which one accedes is that paraconsistent logics are certain sets of formulas (in the present work,~$\seeone$-based theories), axiomatized in certain ways; one does not take the further step that these sets of formulas have any bearing on acceptable reasoning, say, for mathematics or science.  The work presented here cannot refute this view.  Nonetheless, we hope to have provided another weapon in the conceptual arsenal for the paraconsistent researcher and opened a new method for investigating the foundations of paraconsistency.

\bibliographystyle{plain}
\bibliography{paraconsistent}

\end{document}